\documentclass[11pt]{article}
\usepackage{amsmath,amssymb}
\usepackage{url}
\newcommand\nonu{\nonumber}
\newcommand\dstyle\displaystyle

\newcommand\sLP{\\[\smallskipamount]}
\newcommand\mLP{\\[\medskipamount]}
\newcommand\bPP{\\[\bigskipamount]\indent}

\newcommand\ZZ{\mathbb{Z}}
\newcommand\la\lambda
\newcommand\om\omega

\newcommand\half{\frac12}
\newcommand\thalf{\tfrac12}
\newcommand\iy\infty

\newcommand{\qhyp}[5]{\,\mbox{}_{#1}\phi_{#2}\!\left(
  \genfrac{}{}{0pt}{}{#3}{#4};#5\right)}
\newcommand\LHS{left-hand side}
\newcommand\RHS{right-hand side}

\usepackage[pstarrows]{pict2e}
\usepackage{tikz}
\newcommand{\tikzcircle}[2][red,fill=red]{\tikz[baseline=-0.5ex]
  \draw[thick,#1,radius=#2] (0,0) circle ;}%
\newcommand{\black}{\tikzcircle[black, fill=black]{2.5pt}}%
\newcommand{\white}{\tikzcircle[black, fill=white]{2.5pt}}%
\newcommand\td{\tilde}
\newcommand\wt{\widetilde}

\numberwithin{equation}{section}
\newtheorem{theorem}{Theorem}[section]
\newtheorem{proposition}[theorem]{Proposition}

\begin{document}
\title{Charting the $q$-Askey scheme. II. The $q$-Zhedanov scheme}
\author{Tom H. Koornwinder}
\date{\emph{Dedicated to Jaap Korevaar on the occasion of his centennial birthday}}
\maketitle
\begin{abstract}
This is the second in a series of papers which intend to explore
conceptual ways of distinguishing between families in the $q$-Askey scheme and
uniform ways of parametrizing the families. For a system of polynomials $p_n(x)$
in the $q$-Askey scheme satisfying $Lp_n=h_np_n$ with $L$ a second order
$q$-difference operator the $q$-Zhedanov algebra is the algebra generated
by operators $L$ and $X$ (multiplication by $x$). It has two relations in
which essentially five coefficients occur. Vanishing of one or more of
the coefficients corresponds to a subfamily or limit family of the Askey--Wilson
polynomials. An arrow from one family to another means that in the latter family
one more coefficient vanishes. This yields the $q$-Zhedanov scheme given in
this paper.

The $q$-hypergeometric expression of $p_n(x)$ can be interpreted as an
expansion of $p_n(x)$ in terms of certain Newton polynomials.
In our previous paper (Contemporary Math.\ 780)
we used Verde-Star's clean parametrization of such
expansions and we obtained a $q$-Verde-Star scheme, where vanishing of one or
more of these parameters corresponds to a subfamily or limit family.
The actions of the operators $L$ and $X$ on the Newton polynomials can be
expressed in terms
of the Verde-Star parameters, and thus the coefficients for the
$q$-Zhedanov algebra can be expressed in terms of these parameters.
There are interesting differences between the $q$-Verde-Star scheme and
the $q$-Zhedanov scheme, which are discussed in the paper.
\end{abstract}
%
%
%
\section{Introduction}
\textsl{Jaap Korevaar has worked on a wide range of topics in analysis, in
particular complex and harmonic analysis and applications to number theory.
As for special functions, he likes in particular the special functions of number
theory, with the zeta function on top \cite{24}, \cite{25}, but he does certainly not
ignore the classical orthogonal polynomials. In \cite{28}
he used Hermite polynomials, and in work with his PhD student Meyers \cite{26},
\cite{27} he used Legendre and ultraspherical polynomials. The paper \cite{26}
was a key reference in \cite{29}, a paper coauthored by Maryna Viazovska who was
a 2022 Fields Medal winner.}

After pioneering work of Hahn $q$-analogues of classical orthogonal polynomials got
during 1975--1985 an increasing interest. This culminated in the introduction of
the Askey--Wilson polynomials \cite{8}, a 4-parameter family of
$q$-hypergeometric orthogonal polynomials which stand on top of the
$q$-Askey scheme \cite[Ch.~14]{3}. In that scheme there are boxes representing
families of $q$-hypergeometric orthogonal polynomials and arrows connecting
them. An arrow stands for a limit transition or specialization and it
diminishes the number of parameters by at least one.

This paper is the second in a series of papers which intend to explore
conceptual ways of distinguishing between families in the scheme and
uniform ways of parametrizing the families.
First of all we do this for the $q$-Askey scheme,
and tentatively later for its $q=1$ limit
(the original Askey scheme \cite[Ch.~9]{3})
and its $q=-1$ limit \cite{15}, \cite{22}. 
The present and the previous paper
concentrate on the case of general $q$.
The previous paper \cite{1}
was inspired by Verde-Star's paper \cite{9} (see also Vinet and Zhedanov
\cite{16} and, for finite systems, Terwilliger \cite{20}, cf.~\cite[Remark 2.3]{1}),
which classifies the families
in the $q$-Askey scheme in combination with a basis of Newton polynomials
in terms of which they have nice expansions. There are uniform parameters
describing this. The vanishing of one or more of the parameters
corresponds to a subfamily.

The present paper classifies the $q$-Zhedanov algebras \cite[\S3]{2}
associated with families in the $q$-Askey scheme. Such an algebra is generated
by the second order $q$-difference operator $L$ for which the polynomials in the
family are eigenfunctions, and by the operator of multiplication by $x$. Denote
these two operators by $K_1,K_2$, respectively. Then
the algebra is described by two relations which expand certain cubic forms
in $K_1,K_2$ (repeated $q$-commutators) as linear combinations of $K_1$, $K_2$
and the identity operator. The vanishing of one or more of the coefficients
in these linear combinations indicates a subfamily in the $q$-Askey scheme.
The prototypical algebra of this type is the Askey--Wilson algebra,
associated with the Askey--Wilson polynomials. See Zhedanov \cite{13}.

In Verde-Star's approach the action of $L$ on the space of polynomials is
conceptually described by its action on the corresponding Newton polynomials,
while multiplication by $x$ immediately yields an explicit action on the
Newton polynomials.
This makes it possible to express the coefficients of the $q$-Zhedanov algebra
in terms of Verde-Star's parameters. From that it can be read off for which
constraints on Verde-Star's parameters certain $q$-Zhedanov coefficients vanish.
Since we already expressed in \cite{1} Verde-Star's parameters in terms of
the parameters of the corresponding family in the $q$-Askey scheme, we can
associate vanishing properties of the $q$-Zhedanov parameters with
families in the $q$-Askey scheme, and thus the $q$-Zhedanov scheme occurs.
One can find these results in Sections \ref{17}, \ref{90} and \ref{89}, which
form the heart of this paper. It should be observed that our work was anticipated,
in the case of finite systems and by a different approach, by Terwilliger and Vidunas
\cite[Theorem 4.5 and 5.3]{19} and Vidunas \cite{21}.
More limited $q$-Zhedanov schemes, not related to $q$-Verde-Star parameters, were
earlier given by Mazzocco \cite{23}, see also \cite{6}.

Filling in this $q$-Zhedanov scheme gave some surprises. Several families
in the $q$-Askey scheme which usually do not get much attention attain an
independent status in the present scheme. We mention symmetric cases of
Askey--Wilson, Al-Salam--Chihara and big $q$-Jacobi polynomials, and also
$q^2$ versions
of the continuous $q$-Jacobi polynomials and its subfamilies, which give
interesting $q$-polynomials by quadratic transformation. In particular,
continuous $q^2$-Hermite polynomials get after quadratic transformation a
$q$-hypergeometric expression which allows to describe this family by
Verde-Star's approach, something which was not possible in \cite{9} for the
continuous $q$-Hermite polynomials.

Very relevant for this paper is the fact that the Askey--Wilson algebra in its
dependence on the Askey--Wilson parameters $a,b,c,d$ is not only symmetric
in these parameters, but also remains invariant if two parameters $e,f$
from $a,b,c,d$ are sent to $q/e$, $q/f$, respectively. Such a symmetry
under $(a,b,c,d)\to(a,b,qd^{-1},qc^{-1})$
for the Askey--Wilson DAHA \cite{18}, \cite{7}, \cite[\S5.1]{6} was
communicated to me by Marta Mazzocco. She pointed for this to the DAHA symmetry
$(a,b,c,d)\to(a,b,qc^{-1},d)$ observed by Oblomkov \cite[(2.11)]{14}.
Since Askey--Wilson polynomials are not
invariant under these transformations (however, see Section \ref{95}),
it means that a possible charting of
the $q$-Zhedanov scheme cannot be a charting
of the $q$-Askey scheme which completely distinguishes between individual
systems of polynomials.

Duality is also an important theme in this paper, as it was also in the author's
earlier paper~\cite{6} with Mazzocco. From the graphical $q$-Zhedanov scheme in
Section \ref{89} the dualities between the occurring families can be nicely
seen.

It should be emphasized that, although the families from the $q$-Askey scheme
are first of all known as orthogonal polynomials, we completely disregard
orthogonality aspects in this paper. Only those algebraic aspects which give
rise to a $q$-Verde-Star description and a nice associated $q$-Zhedanov algebra
matter here for us.

The contents of this paper are as follows. In Section \ref{93}
we treat some (but not all) families from the $q$-Askey scheme. We concentrate
on the families  which play a particularly important role later in the paper.
Duality properties are already discussed here. This section may be skipped on
first reading. As already said, Section \ref{17} on Verde Star's framework,
Section \ref{90} on using this framework in the context of the $q$-Zhedanov
algebra, and Section \ref{89} which presents the $q$-Zhedanov scheme in
Figure 1, form the heart of this paper. Some subsections of Section \ref{90}
consider examples of $q$-Zhedanov algebras in connection with
$q$-Verde-Star parameters. Here the further symmetries just observed also
enter. But these symmetries appear already in Subsection \ref{59} in connection with
duality, although not yet in connection with the Askey--Wilson algebra.
In Section \ref{88} the duals of continuous $q^2$-Hermite polynomials and
of discrete $q$-Hermite I polynomials are discussed. Finally Section \ref{95}
gives some recapitulation, open questions and thoughts about an optimal
$q$-Askey scheme.
\paragraph{Note}
For definition and notation of $q$-shifted factorials and
$q$-hypergeometric series we follow \cite[\S1.2]{4}.
We will only need terminating series:
\begin{equation*}
\qhyp rs{q^{-n},a_2,\ldots,a_r}{b_1,\ldots,b_s}{q,z}:=
\sum_{k=0}^n \frac{(q^{-n};q)_k}{(q;q)_k}\,
\frac{(a_2,\ldots,a_r;q)_k}{(b_1,\ldots,b_s;q)_k}\,
\big((-1)^k q^{\half k(k-1)}\big)^{s-r+1}z^k.
\end{equation*}
Here
$(b_1,\ldots,b_s;q)_k:=(b_1;q)_k\ldots(b_s;q)_k$ with
$(b;q)_k:=(1-b)(1-qb)\ldots(1-q^{k-1}b)$ the
$q$-shifted factorial.
\section{Preliminaries about polynomials in the $q$-Askey scheme}
\label{93}
\subsection{Askey--Wilson polynomials and some subfamilies}
Assume $a\ne0$ and $ab,ac,ad,abcd\ne1,q^{-1},q^{-2},\ldots$\;.
Define \emph{Askey--Wilson polynomials} by
\begin{equation}
R_n(z;a,b,c,d\,|\,q)=R_n(z):=
\qhyp43{q^{-n},q^{n-1}abcd,az,az^{-1}}{ab,ac,ad}{q,q}.
\label{27}
\end{equation}
Then $R_n(z)$ is a symmetric Laurent polynomial in $z$ of degree $n$.
In usual notation and normalization \cite[(14.1.1)]{3}
the Askey--Wilson polynomials are written as
\[
p_n\big(\thalf(z+z^{-1});a,b,c,d\,|\,q\big):=
a^{-n}(ab,ac,ad;q)_n\,R_n(z;a,b,c,d\,|\,q).
\]
In this form they are symmetric in $a,b,c,d$ (see \cite[p.6]{8}).
In the normalization \eqref{27} the polynomials are still symmetric in
$b,c,d$, but for $a\leftrightarrow b$ we have
\begin{equation}
R_n(z;a,b,c,d\,|\,q)=\left(\frac ab\right)^n \frac{(bc,bd;q)_n}{(ac,ad;q)_n}\,
R_n(z;b,a,c,d\,|\,q),
\label{72}
\end{equation}
where we assume in addition that $b\ne0$ and $bc,bd\ne1,q^{-1},q^{-2},\ldots$\;.
For $c=-a$, $d=-b$ we have the
\emph{symmetric Askey--Wilson polynomials} \cite[p.7]{8} with the symmetry
\begin{equation}
R_n(-z;a,b,-a,-b\,|\,q)=(-1)^n R_n(z;a,b,-a,-b\,|\,q).
\label{28}
\end{equation}

The polynomials $R_n(z)$ given by \eqref{27} satisfy the eigenvalue equation
\cite[(14.1.7)]{3}
\begin{equation}
(LR_n)(z)=(q^{-n}+abcd q^{n-1})R_n(z),
\label{54}
\end{equation}
where
\begin{multline}
(Lf)(z)=
\big(1+q^{-1}abcd\big)f(z)+\frac{(1-az)(1-bz)(1-cz)(1-dz)}{(1-z^2)(1-qz^2)}\,
\bigl(f(qz)-f(z)\bigr)\\
+\frac{(a-z)(b-z)(c-z)(d-z)}{(1-z^2)(q-z^2)}\,
\bigl(f(q^{-1}z)-f(z)\bigr).
\label{55}
\end{multline}
They also satisfy the three-term recurrence relation
(for $n=1,2,\ldots$)
\begin{multline}
(z+z^{-1})R_n(z)=(a+a^{-1})R_n(z)\\*
+\frac{(1-q^nab)(1-q^nac)(1-q^nad)(1-q^{n-1}abcd)}
{a(1-q^{2n-1}abcd)(1-q^{2n}abcd)}\,(R_{n+1}(z)-R_n(z))\\*
+\frac{a(1-q^n)(1-q^{n-1}bc)(1-q^{n-1}bd)(1-q^{n-1}cd)}
{(1-q^{2n-2}abcd)(1-q^{2n-1}abcd)}\,(R_{n-1}(z)-R_n(z)).
\label{65}
\end{multline}
The starting values $R_0(z)=1$ and $R_1(z)$ as given by \eqref{27} satisfy
\eqref{65} in its reduced form
\[
(z+z^{-1})R_0(z)=(a+a^{-1})R_0(z)\\
+\frac{(1-ab)(1-ac)(1-ad)}
{a(1-abcd)}\,(R_1(z)-R_0(z)).
\]
Note that here the factors $(1-q^{-1}abcd)$ in numerator and
denominator in \eqref{65} cancel, which even turns out to be correct if $abcd=q$.

In the following we mention a few special cases
of Askey--Wilson which will play a remarkable role in the $q$-Zhedanov scheme.
These special cases have their own, generally used, notations
(see \cite[Ch.~14]{3}), but we will just denote them as a restricted $R_n(z)$.
\paragraph{Continuous dual $q$-Hahn}
\cite[(14.3.1)]{3}
\begin{equation}
R_n(z;a,b,c,0\,|\,q)=\qhyp32{q^{-n},az,az^{-1}}{ab,ac}{q,q}.
\label{83}
\end{equation}
\paragraph{Al-Salam--Chihara}
\cite[(14.8.1)]{3}
\begin{equation}
R_n(z;a,b,0,0\,|\,q)=\qhyp32{q^{-n},az,az^{-1}}{ab,0}{q,q}.
\label{82}
\end{equation}
\paragraph{Symmetric Al-Salam--Chihara}
\begin{equation}
R_n(z;a,-a,0,0\,|\,q)=\qhyp32{q^{-n},az,az^{-1}}{-a^2,0}{q,q}.
\label{85}
\end{equation}
For these polynomials the symmetry \eqref{28} holds with $b=0$.
\paragraph{Continuous $q$-Jacobi polynomials}
\cite[(14.10.1)]{3}
\begin{equation}
R_n(z;a,-b,q^\half a,-q^\half b\,|\,q)=
\qhyp43{q^{-n},q^n a^2b^2,az,az^{-1}}{-ab,q^\half a^2,-q^\half ab}{q,q}.
\label{29}
\end{equation}
\paragraph{Continuous $q$-Laguerre polynomials}
\cite[(14.18.1)]{3},
\begin{equation}
R_n(z;a,q^\half a,0,0\,|\,q)
=\qhyp32{q^{-n},az,az^{-1}}{q^\half a^2,0}{q,q}.
\label{30}
\end{equation}
\paragraph{Continuous $q$-ultraspherical polynomials}
\cite[(14.10.17)]{3}
\begin{equation}
R_n(z;a,-a,q^\half a,-q^\half a\,|\,q)
=\qhyp43{q^{-n},q^n a^4,az,az^{-1}}{-a^2,q^\half a^2,-q^\half a^2}{q,q}.
\label{31}
\end{equation}
For these polynomials the symmetry \eqref{28} holds with $b=q^\half a$.
\paragraph{Continuous $q$-Hermite polynomials}
\cite[\S14.26]{3},
\begin{equation}
H_n\big(\thalf(z+z^{-1})\,|\,q\big):=\begin{cases}
\dstyle\lim_{a\to0} a^{-n} R_n(z;a,q^\half a,0,0\,|\,q),\mLP
\dstyle\lim_{a\to0} a^{-n} R_n(z;a,-a,q^\half a,-q^\half a\,|\,q).
\end{cases}
\label{32}
\end{equation}
This is no longer a special case of Askey--Wilson, but a limit case.
There is the symmetry $H_n(-x)=(-1)^n H_n(x)$.
\paragraph{$q\to q^2$ transformations}
The transformation \cite[(3.10.13)]{4}
\begin{equation*}
\qhyp43{a^2,b^2,c,d}{q^\half ab,-q^\half ab,-cd}{q,q}
=\qhyp43{a^2,b^2,c^2,d^2}{qa^2b^2,-cd,-qcd}{q^2,q^2},
\end{equation*}
together with \eqref{28}, allows to express the polynomials
\eqref{29}--\eqref{32} for base $q^2$ in a different way:
\begin{align}
R_n(z;a,-b,q^\half,-q^\half\,|\,q)&=R_n(z;a,-b,qa,-qb\,|\,q^2),
\label{33}\\*
q^{-\half n}\,\frac{(-q;q)_n}{(-q^\half a;q)_n}\,
R_n(z;q^\half,-q^\half,a,0\,|\,q)&=a^{-n} R_n(z;a,qa,0,0\,|\,q^2),
\label{34}\\*
q^{-\half n}\,\frac{(-q;q)_n}{(-a^2;q)_n}\,
R_n(z;q^\half,-q^\half,a,-a\,|\,q)&=a^{-n}R_n(z;a,-a,qa,-qa\,|\,q^2),
\label{35}\\*
q^{-\half n}(-q;q)_n\,R_n(z;q^\half,-q^\half,0,0\,|\,q)&=
H_n\big(\thalf(z+z^{-1})\,|\,q^2\big).
\label{36}
\end{align}
We call the polynomials $R_n(z)$ on the left-hand sides
the continuous $q^2$-Jacobi, continuous $q^2$-Laguerre,
continuous $q^2$-ultraspherical and continuous $q^2$-Hermite
polynomials, respectively. The continuous $q^2$-Jacobi polynomials as
represented by the \LHS\ of \eqref{33} were first observed by M.~Rahman, see
\cite[p.468]{3}.
\subsection{Big $q$-Jacobi polynomials and some subfamilies}
Assume $ab,a,c\ne q^{-1},q^{-2},\ldots$\;.
Define \emph{big $q$-Jacobi polynomials} \cite[(14.5.1)]{3} by
\begin{equation}
P_n(x;a,b,c;q)=P_n(x):=\qhyp32{q^{-n},q^{n+1}ab,x}{qa,qc}{q,q}.
\label{48}
\end{equation}
They are limit cases \cite[(14.1.18)]{3} of Askey--Wilson polynomials:
\begin{equation}
P_n(x;a,b,c;q)=\lim_{\la\to0}
R_n(\la^{-1}x;\la,\la bc^{-1},q\la^{-1}a,q\la^{-1}c\,|\,q).
\label{73}
\end{equation}
The $S_4$ symmetry of Askey--Wilson in its parameters reduces to an
 $S_2\times S_2$
symmetry for big $q$-Jacobi. Here the $S_3$ symmetry which is obvious from
\eqref{27} reduces to a symmetry obvious from \eqref{48}:
\begin{equation}
\begin{split}
P_n(x;a,b,c;q)&=P_n(x;c,abc^{-1},a;q)\quad\mbox{or}
\\
P_n(x;a,bc,c;q)&=P_n(x;c,ba,a;q)\quad\mbox{(equivalently)}.
\end{split}
\label{67}
\end{equation}
The other symmetry is a limit case of \eqref{67} under the limit \eqref{73}:
\begin{equation}
\begin{split}
&P_n(x;a,b,c;q)=\left(\frac cb\right)^n
\frac{(qabc^{-1},qb;q)_n}{(qa,qc;q)_n}\,
P_n(bc^{-1}x;abc^{-1},c,b;q)\quad\mbox{or}\\
&c^{-n}(qac,qc;q)_n P_n(cx;ac,b,c;q)=b^{-n}(qab,qb;q)_n P_n(bx;ab,c,b;q)\quad
\mbox{(equivalenty)},
\end{split}
\label{68}
\end{equation}
where we assume that $b,c\ne0$ and (in the first equality)
$a,b,c,abc^{-1}\ne q^{-1},q^{-2},\ldots$\;.
In particular,  we have the \emph{symmetric big $q$-Jacobi polynomials}
\begin{equation}
P_n(x;a,a,-a;q)=\qhyp32{q^{-n},q^{n+1}a^2,x}{qa,-qa}{q,q},
\label{86}
\end{equation}
which satisfy
\[
P_n(-x;a,a,-a;q)=(-1)^n P_n(x;a,a,-a;q)
\]
by \eqref{68} and \eqref{67}.

The polynomials $P_n(x)$ given by \eqref{48} satisfy the eigenvalue equation
\cite[(14.5.5)]{3}
\begin{equation}
(LP_n)(x)=(q^{-n}+abq^{n+1})P_n(x),
\end{equation}
where
\begin{equation}
(Lf)(x)=
(1+qab)f(x)+qa(x-1)(bx-c)\big(f(qx)-f(x)\big)
+(x-qa)(x-qc)\big(f(q^{-1}x)-f(z)\big).
\end{equation}
They also satisfy the three-term recurrence relation
(for $n=1,2,\ldots$)
\begin{multline}
xP_n(x)=P_n(x)+\frac{(1-q^{n+1}a)(1-q^{n+1}ab)(1-q^{n+1}c)}
{(1-q^{2n+1}ab)(1-q^{2n+2}ab)}\,(P_{n+1}(x)-P_n(x))\\
-q^{n+1}ac\,\frac{(1-q^n)(1-q^nabc^{-1})(1-q^nb)}
{(1-q^{2n}ab)(1-q^{2n+1}ab)}\,(P_{n-1}(x)-P_n(x)).
\label{69}
\end{multline}
The starting values $P_0(x)=1$ and $P_1(x)$ as given by \eqref{48} satisfy
\eqref{69} in its reduced form
\[
xP_0(x)=P_0(x)+\frac{(1-qa)(1-qc)}{1-q^2ab}\,(P_1(x)-P_0(x)).
\]

We mention some further subfamilies of the big $q$-Jacobi polynomials.
\paragraph{Little $q$-Jacobi}
\cite[(14.12.1) and p.442, Remarks]{3}\\
These can be seen as a special case of big
$q$-Jacobi polynomials:
\begin{equation}
P_n(x;a,b,0;q)=
\qhyp32{q^{-n},q^{n+1}ab,x}{qa,0}{q,q}.
\label{70}
\end{equation}
However, they are usually notated and defined as
\begin{equation}
p_n(x;a,b;q):=\qhyp21{q^{-n},q^{n+1}ab}{qa}{q,qx}.
\label{71}
\end{equation}
The expressions in \eqref{70} and \eqref{71} are related by
\begin{equation}
(-a)^n q^{\half n(n+1)}\,\frac{(qb;q)_n}{(qa;q)_n}\,
p_n((qa)^{-1}x;b,a;q)=P_n(x;a,b,0;q).
\label{87}
\end{equation}
\paragraph{Big $q$-Laguerre}
\cite[14.11.1)]{3}
\begin{equation}
P_n(x;a,0,c;q)=\qhyp32{q^{-n},x,0}{qa,qc}{q,q}.
\label{84}
\end{equation}
\subsection{Duality}
\label{59}
In our charting of the $q$-Askey scheme we consider systems of polynomials
$\{p_n(x)\}$ in the scheme up to dilation of $x$ by a nonzero constant.
In particular, taking $-1$ for this constant and noting from \eqref{27} that
\begin{equation}
R_n(-z;a,b,c,d\,|\,q)=R_n(z;-a,-b,-c,-d\,|\,q),
\end{equation}
parameter values $a,b,c,d$ and $-a,-b,-c,-d$ will be identified for
charting purposes.

Define {\em dual parameters}
$\td a,\td b,\td c,\td d$ in terms of the Askey--Wilson parameters
$a,b,c,d$ by
\begin{equation}
\td a=(q^{-1}abcd)^\half,\quad \td b=ab/\td a,\quad \td c=ac/\td a,\quad
\td d=ad/\td a,
\label{53}
\end{equation}
see \cite[\S\S~5.7, 8.5]{7} and \cite[\S2.3]{6}.
Because of the square root in the definition of $\td a$, equation \eqref{53}
defines the values of the dual parameters up to possible common multiplication
by $-1$. So for charting purposes this is harmless.

From \eqref{27} we have the duality relation
\begin{equation}
R_n(a^{-1}q^{-m};a,b,c,d\,|\,q)=
R_m(\td a^{-1} q^{-n};\td a,\td b,\td c,\td d\,|\,q)\qquad
(m,n\in\ZZ_{\ge0}),
\label{37}
\end{equation}
with both sides equal to
\begin{equation}
\qhyp43{q^{-n},q^{n-1}abcd,q^{-m},q^m a^2}{ab,ac,ad}{q,q}.
\label{39}
\end{equation}
Formula \eqref{37} gives a duality between the two Askey--Wilson polynomials
$R_n(z;a,b,c,d\,|\,q)$ and $R_m(w;\td a,\td b,\td c,\td d\,|\,q)$.
\paragraph{Remark}
The above concept of duality is rather weak. It certainly does not
involve in general that orthogonality relations for the original system
imply orthogonality for the dual system. See \cite{10} for a treatment
with many examples of dual systems in this stronger sense. Still our
weak notion of duality, when occurring in the $q$-Askey scheme, usually
makes it possible to relate the three-term recurrence relation and the
second order $q$-difference eigenvalue equation for the one system
in a formal way with the eigenvalue equation and the recurrence relation,
respectively, for the dual system. For instance, the recurrence relation
\eqref{65} for Askey--Wilson polynomials can be rewritten as
\begin{multline}
a(z+z^{-1})R_n(z)=(1+q^{-1}\td a\td b\td c\td d)R_n(z)\\*
+\frac{(1-q^n\td a^2)(1-q^n\td a\td b)(1-q^n\td a\td c)(1-q^n\td a\td d)}
{(1-q^{2n}\td a^2)(1-q^{2n+1}\td a^2)}\,(R_{n+1}(z)-R_n(z))\\*
+\frac{(\td a-q^n\td a)(\td b-q^n\td a)(\td c-q^n\td a)(\td d-q^n\td a)}
{(1-q^{2n}\td a^2)(q-q^{2n}\td a^2)}\,(R_{n-1}(z)-R_n(z)).
\end{multline}
Compare this with \eqref{55} with $z$ replaced by $q^n\td a$ in the coefficients
on the \RHS.
\bPP
Note that the definition of the dual parameters depends on a choice of one of
the original parameters. In \eqref{53} this is $a$. The duals of the dual
parameters $\td a, \td b, \td c, \td d$ with respect to $\td a$ are
the original
parameters $a,b,c,d$ again, up to common multiplication by $-1$.
(In the following this identification of parameters with their opposites
will be silently assumed.)\;
Now let $a',b',c',d'$ be dual parameters of
$a,b,c,d$ with respect to $b$:
\begin{equation*}
a'=\frac{ab}{(q^{-1}abcd)^\half}\,,\quad
b'=(q^{-1}abcd)^\half,\quad
c'=\frac{bc}{(q^{-1}abcd)^\half}\,,\quad
d'=\frac{bd}{(q^{-1}abcd)^\half}\,.
\end{equation*}
Then the dual parameters of $a',b',c',d'$ with respect to $a'$ are
$b,a,qd^{-1},qc^{-1}$. This result can be alternatively stated by
the commutative diagram
\[
\begin{matrix}
a\;b\;c\;d&\xrightarrow{\sim}&\td a\;\td b\;\td c\;\td d\\
\downarrow&&\downarrow\\
a\;b\;\frac qd\;\frac qc&\xrightarrow{\sim}&\td b\;\td a\;\td c\;\td d
\end{matrix}
\]
From this we have:
\begin{proposition}\label{62}
Let $f(a,b,c,d)$ be a symmetric function which is also invariant under common
multiplication of $a,b,c,d$ by $-1$.
Then $f$ is symmetric in
$\td a,\td b,\td c,\td d$ iff $f$ is invariant under any mapping which sends
an even number of its variables $a,b,c,d$ to $q/a,q/b,q/c,q/d$, respectively.
\end{proposition}
The Askey--Wilson polynomials, although invariant under permutations of
$a,b,c,d$, are no longer invariant under the larger group which also involves
mappings sending two parameters $e,f$ out of $a,b,c,d$
to $qe^{-1},qf^{-1}$, respectively, but such symmetries
will play a role later
in this paper, see Subsection \ref{60}.
The larger group is isomorphic to the semidirect
product $S_4\ltimes (\ZZ/2\ZZ)^3$
(where $(\ZZ/2\ZZ)^3$ is the group of an even number of sign changes of
$1,2,3,4$),
which is also the Weyl group for the root system $\textup{D}_4$, see
\cite[\S12.1]{17}.

Below we list some dualities obtained from \eqref{37}, \eqref{39}
by specialization or (possibly scaled) limit.
\paragraph{Continuous dual $q$-Hahn $\longleftrightarrow$ big $q$-Jacobi}
\cite[(44)]{6}
\begin{align}
&R_n(z;a,b,c,0\,|\,q)\;\longleftrightarrow\;
P_m(y;q^{-1}ab,ab^{-1},q^{-1}ac;q)\label{74}\\
&R_n(a^{-1}q^{-m};a,b,c,0\,|\,q)=P_m(q^{-n};q^{-1}ab,ab^{-1},q^{-1}ac;q)=
\qhyp32{q^{-n},q^{-m},q^ma^2}{ab,ac}{q,q}.\nonu
\end{align}
Note that there is a one-to-one correspondence
$\pm(a,b,c)\leftrightarrow (q^{-1}ab,ab^{-1},q^{-1}ac)$.
By \eqref{72} there is a symmetry $(a,b,c)\leftrightarrow(b,a,c)$ for
suitably normalized continuous dual $q$-Hahn.
By the duality \eqref{74} this leads to a symmetry
$(a,b,c)\leftrightarrow(a,b^{-1},b^{-1}c)$ for the big $q$-Jacobi parameters
$a,b,c$, although it does not leave $P_n(x;a,b,c;q)$ invariant.
Conversely, the symmetry \eqref{68} for big $q$-Jacobi polynomials leads to
a symmetry $(a,b,c)\leftrightarrow (a,qc^{-1},qb^{-1})$ for the continuous
dual Hahn parameters $a,b,c$, which, however, does not leave
$R_n(z;a,b,c,0\,|\,q)$ invariant.
\paragraph{Al-Salam--Chihara $\longleftrightarrow$ little $q$-Jacobi}
\cite[(75)]{6}
\begin{align}
&R_n(z;a,b,0,0\,|\,q)\;\longleftrightarrow\;
P_m(y;q^{-1}ab,ab^{-1},0;q)\\
&R_n(a^{-1}q^{-m};a,b,0,0\,|\,q)=P_m(q^{-n};q^{-1}ab,ab^{-1},0;q)=
\qhyp32{q^{-n},q^{-m},q^ma^2}{ab,0}{q,q}.\nonu
\end{align}
\paragraph{Continuous $q^2$-Jacobi $\longleftrightarrow$ symmetric Askey--Wilson}
\begin{align}
&R_n(z;a,-b,q^\half,-q^\half\,|\,q)\;\longleftrightarrow\;
R_m(w;(ab)^\half,-(ab)^\half,(qab^{-1})^\half,-(qab^{-1})^\half\,|\,q),
\label{38}\\
&R_n(a^{-1}q^{-m};a,-b,q^\half,-q^\half\,|\,q)=
R_m((ab)^{-\half}q^{-n};(ab)^\half,-(ab)^\half,(qab^{-1})^\half,
-(qab^{-1})^\half\,|\,q)\nonu\\*
&\qquad\qquad\qquad\qquad\qquad\qquad\qquad\qquad
=\qhyp43{q^{-n},q^n ab,q^{-m},q^m a^2}{-ab,q^\half a,-q^\half a}{q,q}.
\nonu
\end{align}
\paragraph{Continuous $q^2$-Laguerre $\longleftrightarrow$ symmetric
big $q$-Jacobi}
\begin{align}
&R_n(z;a,q^\half,-q^\half,0\,|\,q)\;\longleftrightarrow\;
P_m(y;q^{-\half}a,q^{-\half}a,-q^{-\half}a;q),
&\label{40}\\
&R_n(a^{-1}q^{-m};a,q^\half,-q^\half,0\,|\,q)=
P_m(a^{-1}q^{-n-\half};q^{-\half}a,q^{-\half}a,1,1;q)
=\qhyp32{q^{-n},q^{-m},q^m a^2}{q^\half a,-q^\half a}{q,q}.\nonu
\end{align}
\paragraph{Continuous $q^2$-ultraspherical self-dual}
\begin{align}
&R_n(z;a,-a,q^\half,-q^\half\,|\,q)\;\longleftrightarrow\;
R_m(w;a,-a,q^\half,-q^\half\,|\,q),
\label{41}\\
&R_n(a^{-1}q^{-m};a,-a,q^\half,-q^\half\,|\,q)=
R_m(a^{-1}q^{-n};a,-a,q^\half,-q^\half\,|\,q)\nonu\\
&\qquad\qquad\qquad\qquad\qquad\qquad\qquad\qquad
=\qhyp43{q^{-n},q^n a^2,q^{-m},q^m a^2}{-a^2,q^\half a,-q^\half a}{q,q}.
\nonu
\end{align}
\section{Verde-Star's description of the $q$-Askey scheme}\label{17}
Monic polynomials $u_n$ in the $q$-Askey scheme can be described by the data
\cite{9}, \cite[(3.1)--(3.4)]{1}:
\begin{align}
&u_n(x)=\sum_{k=0}^n c_{n,k}\,v_k(x),\quad
v_k(x)=\prod_{j=0}^{k-1}(x-x_j),\quad
c_{n,k}=\prod_{j=k}^{n-1}\frac{g_{j+1}}{h_n-h_j}\,,\label{1}\\
&x_k=b_1 q^k+b_2 q^{-k},\quad
h_k=a_1 q^k+a_2 q^{-k},\quad
g_k=d_3 q^{2k}+d_1 q^k+d_0+d_2 q^{-k}+d_4 q^{-2k},\label{2}\\
&\sum_{i=0}^4 d_i=0,\quad d_3=q^{-1}a_1b_1,\quad d_4=qa_2b_2,\label{3}\\
&a_2\ne a_1 q^{\ZZ_{>0}},\quad\mbox{in particular,\quad $a_1$ or $a_2\ne 0$},
\qquad \mbox{$d_i\ne0$ for some $i$}.\label{4}
\end{align}

Define a linear operator $L$ on the space of polynomials by
\begin{equation}
Lv_0=0,\qquad Lv_n=h_nv_n+g_nv_{n-1},\quad n>0,
\label{5}
\end{equation}
or equivalently, in view of \eqref{1}, by
\begin{equation}
Lu_n=h_nu_n,\quad n\ge0.
\label{6}
\end{equation}
For specific $u_n$ in the $q$-Askey scheme the eigenvalue equation
\eqref{5}, \eqref{6}
can be rewritten as a second order $q$-difference equation (the
generalized Bochner property \cite{5}).

In \cite[(3.2)]{1} $x_k$ had an additional term $b_0$ and $h_k$ had an additional
term $a_0$, but without loss of generality these can be omitted, because
the term $b_0$ only leads to a translation of the $x$-variable and the term $a_0$
only leads to adding a constant to the operator $L$ in \eqref{5}.

The polynomials $u_n$ are determined by the 9 parameters
$a_1,a_2,b_1,b_2,d_0,d_1,d_2,d_3,d_4$ under the 3 constraints \eqref{3}.
There are two other operations on these parameters which do not lead to
essential changes.
Simultaneous multiplication of $a_1,a_2$ and $d_0,d_1,d_2,d_3,d_4$ by $\mu\ne0$
only leads to multiplication of $L$ by $\mu$.
Simultaneous multiplication of 
$b_1,b_2$ and $d_0,d_1,d_2,d_3,d_4$ by $\rho\ne0$ leads to replacing $u_n(x)$
by $\rho^n u_n(\rho^{-1}x)$.
So essentially there are only 4 parameters, just as the Askey--Wilson
polynomials (corresponding to the generic case of these parameters) have
4 parameters.

As observed in \cite[\S3]{1},
there are two further remarkable operations which can be performed on the
11 parameters:
\mLP
\emph{$q\leftrightarrow q^{-1}$ exchange}:\quad
$a_1\leftrightarrow a_2$,\;
$b_1\leftrightarrow b_2$,\;
$d_1\leftrightarrow d_2$,\;
$d_3\leftrightarrow d_4$.
\mLP
\emph{$x\leftrightarrow h$ duality}:\quad
$a_1\leftrightarrow b_1$,\;
$a_2\leftrightarrow b_2$; assume also that
$b_2\ne b_1 q^{\ZZ_{>0}}$, in particular, $b_1$ or $b_2\ne0$.
This relates $u_n$ given by \eqref{1} to its \emph{dual} $\wt u_n$
given by
\begin{equation}
\wt u_n(x)=\sum_{k=0}^n \wt c_{n,k}\,\wt v_k(x),\quad
\wt v_k(x)=\prod_{j=0}^{k-1}(x-h_j),\quad
\wt c_{n,k}=\prod_{j=k}^{n-1}\frac{g_{j+1}}{x_n-x_j}\,.\label{16}
\end{equation}
If we put
\begin{align}
U_n(x)&:=\prod_{j=0}^{n-1} \frac{h_n-h_j}{g_{j+1}}\times u_n(x)
\;\;=\sum_{k=0}^n\prod_{j=0}^{k-1} \frac{(h_n-h_j)(x-x_j)}{g_{j+1}}\,,\label{50}\\*
\wt U_m(y)&:=\prod_{j=0}^{m-1} \frac{x_m-x_j}{g_{j+1}}\times \wt u_m(y)
=\sum_{k=0}^m\prod_{j=0}^{k-1} \frac{(x_m-x_j)(y-h_j)}{g_{j+1}}\,,
\label{51}
\end{align}
then
\begin{equation}
U_n(x_m)=\wt U_m(h_n).
\label{52}
\end{equation}
If $U_n$ or $\tilde U_m$ are written in a special case with \eqref{2}
substituted then they will appear as $q$-hypergeometric series.
\section{$q$-Zhedanov algebra with $q$-Verde-Star parameters}
\label{90}
For $\{u_n\}$ a system of polynomials in the $q$-Askey scheme consider operators
$K_1$, $K_2$ acting on the space of polynomials in one variable $x$,
$K_1$ being the second order $q$-difference operator which has the $u_n$ as
eigenfunctions and $K_2$ the operator of multiplication by $x$. Then
\begin{equation}\label{12}
\begin{split}
&(q+q^{-1})K_2K_1K_2-K_2^2K_1-K_1K_2^2=C_1K_1+DK_2+G_1,\\
&(q+q^{-1})K_1K_2K_1-K_1^2K_2-K_2K_1^2=C_2K_2+DK_1+G_2,
\end{split}
\end{equation}
for certain constants $C_1,C_2,D,G_1,G_2$.
This is a rewriting of \cite[(3.2)]{2} by eliminating there $K_3$ and
by substituting $R=1-\half(q+q^{-1})$.
From \cite[(3.2)]{2} there would have been additional terms $A_2K_2^2$
and $A_1K_1^2$ on the right-hand sides of the two equalities in \eqref{12},
respectively. But these can be removed by adding suitable constants to $K_1$
and $K_2$.
There is also a Casimir operator $Q$, commuting with $K_1$ and $K_2$
and given by \cite[(3.4)]{2}:
\begin{multline}\label{13}
Q=-\half(q+q^{-1})(K_1K_2^2K_1+K_2K_1^2K_2)+K_1K_2K_1K_2+K_2K_1K_2K_1\\*
+\half(q+q^{-1})(C_1K_1^2+C_2K_2^2)+D(K_1K_2+K_2K_1)
+\half(2+q+q^{-1})(G_1K_1+G_2K_2),
\end{multline}
Then, for a certain constant $\om$,
\begin{equation}\label{15}
Q=\om.
\end{equation}
Note the \emph{duality} symmetry for \eqref{12}, \eqref{13}:
\begin{equation}
K_1\leftrightarrow K_2,\quad
C_1\leftrightarrow C_2,\quad
G_1\leftrightarrow G_2.
\label{19}
\end{equation}

The operators $K_1$ and $K_2$ act on the Newton polynomials $v_n$
(see \eqref{1} and \eqref{5}) by
\[
K_1v_n=h_nv_n+g_nv_{n-1},\qquad
K_2v_n=v_{n+1}+x_nv_n.
\]
So by the correspondence $f=\{f_n\}_{n=0}^\iy\longleftrightarrow
\sum_n f_n v_n(x)$ these operators also act on sequences $f$ by
\begin{equation}
(K_1 f)_n=h_n f_n+g_{n+1} f_{n+1},\qquad (K_2 f)_n=x_n f_n+f_{n-1}.
\label{10}
\end{equation}
Substitution of \eqref{10} on the left-hand sides of \eqref{12}
leads to the right-hand sides with the constants explicitly expressed in terms
of the parameters $a_1,a_2,b_1,b_2,d_0,d_1,d_2$. We get 
\begin{equation}\label{14}
\begin{split}
C_1&=(q-q^{-1})^2b_1b_2,\qquad
C_2=(q-q^{-1})^2a_1a_2,\\
D\;&=-(q^\half-q^{-\half})^2\,
\big((q^{-\half}a_1-q^\half a_2)(q^{-\half} b_1-q^\half b_2)+d_1+d_2\big),\\
G_1&=(q^\half-q^{-\half})(q-q^{-1})(q^{-\half}b_1d_2+q^\half b_2d_1),\\
G_2&=(q^\half-q^{-\half})(q-q^{-1})(q^{-\half}a_1d_2+q^\half a_2d_1),
\end{split}
\end{equation}
We can also express $\om$ in \eqref{15} in terms of these parameters:
\begin{align}
\om=(q+q^{-1}-2)\Big(\big(q^{-1}a_1^2+(q+q^{-1})a_1a_2+qa_2^2\big)
\big(q^{-1}b_1^2+(q+q^{-1})b_1b_2+qb_2^2\big)\nonu\\*
+\big(d_1^2+d_2^2-(q+q^{-1})d_1d_2\big)
+(q+q^{-1})(d_1+d_2)(a_1b_2+a_2b_1)\nonu\\*
+2(d_1+d_2)(q^{-1}a_1b_1+qa_2b_2)\Big).\label{18}
\end{align}

Since multiplication of $K_1$ by a nonzero constant $\mu$ and multiplication
of $K_2$ by a nonzero constant $\rho$ does not essentially change the
algebra generated by $K_1$ and $K_2$, it follows from \eqref{12}, \eqref{13} and
\eqref{15} that this algebra is essentially determined
by
\begin{equation}
\{\rho^2 C_1,\,\mu^2C_2,\,\rho\mu D,\,\rho^2\mu\,G_1,\,
\rho\mu^2\,G_2,\,\rho^2\mu^2 \om\}
\label{91}
\end{equation}
with $\mu,\rho$ arbitrarily nonzero.
By \eqref{14}, \eqref{18} this precisely matches with the allowed multiplication
by constants of the $q$-Verde-Star parameters in Section \ref{17}.
So again there are essentially four parameters coming from the $q$-Zhedanov
coefficients,
just as the Askey--Wilson polynomials have four parameters.

The constants $C_1$, $C_2$ and $D$ are invariant under
$q\leftrightarrow q^{-1}$,
but not so $G_1$, $G_2$ and $\om$. The property for $G_1$ or $G_2$ to be
zero for all $q$ is clearly invariant under  $q\leftrightarrow q^{-1}$.
If $q\leftrightarrow q^{-1}$ is combined with
$a_1\leftrightarrow a_2$,\;
$b_1\leftrightarrow b_2$,\;
$d_1\leftrightarrow d_2$,\;
$d_3\leftrightarrow d_4$,
then all constants in \eqref{14}, \eqref{18} remain unchanged.

By \eqref{14}, \eqref{18} the duality symmetry 
$K_1\leftrightarrow K_2$, $C_1\leftrightarrow C_2$, $G_1\leftrightarrow G_2$
in \eqref{19}
is compatible with the $x\leftrightarrow h$ duality
$a_1\leftrightarrow b_1$, $a_2\leftrightarrow b_2$ in Section \ref{17}.

In the following subsections we consider some special cases of the $q$-Zhedanov
algebra in connection with the $q$-Verde-Star parameters: the general case
associated with the Askey--Wilson polynomials, followed by a discussion of
the $(c,d)\to(qd^{-1},qc^{-1})$ symmetry, and finally the dual cases associated
with the continuous dual $q$-Hahn and the big $q$-Jacobi polynomials.
\subsection{The Askey--Wilson algebra}
The generic case of the $q$-Zhedanov algebra is the algebra associated with the
Askey--Wilson polynomials. Put
$U_n(z+z^{-1})=R_n(z;a,b,c,d\,|\,q)$. Then this
can be written in the form \eqref{51} with
\begin{equation}
\begin{split}
&x_k=aq^k+a^{-1}q^{-k},\qquad
h_k=q^{-k}+abcdq^{k-1},\\
&g_k=q^{-2k+1}a^{-1}(1-abq^{k-1})(1-acq^{k-1})(1-adq^{k-1})(1-q^k),
\end{split}
\label{57}
\end{equation}
so
\begin{equation}
\begin{split}
&b_1=a,\quad
b_2=a^{-1}\quad
a_1=q^{-1}abcd,\quad
a_2=1,\\
&d_1=-q^{-2}a(abcd+q(bc+bd+cd)),\quad
d_2=-(b+c+d+qa^{-1}).
\end{split}
\label{61}
\end{equation}
Let $e_1,e_2,e_3,e_4$ be the elementary symmetric polynomials in
$a,b,c,d$:
\begin{equation}
\begin{split}
&e_1=a+b+c+d,\qquad
e_2=ab+ac+bc+ad+bd+cd,\\
&\qquad\quad
e_3=abc+abd+acd+bcd,\qquad
e_4=abcd.
\end{split}
\end{equation}
Then, by \eqref{14} and \eqref{18},
\begin{equation}
\begin{split}
&C_1 =(q-q^{-1})^2,\quad
C_2 =q^{-1}(q-q^{-1})^2 e_4,\quad
D =(1-q^{-1})^2(e_3+qe_1),\\
&G_1 =-q^{-3}(1-q)^2(1+q)(e_4+qe_2+q^2),\quad
G_2 =-q^{-3}(1-q)^2(1+q)(e_1e_4+qe_3),
\end{split}
\label{56}
\end{equation}
and
\begin{multline}
\om=(q+q^{-1}-2)\Big(e_1^2+q^{-2}e_3^2-(1+q^{-2})e_1e_3-q^{-3}(1+q)^2 e_2e_4\\
+q^{-3}(1-q^2)^2e_4-q^{-1}(1+q)^2e_2\Big).
\label{58}
\end{multline}
In \cite[(16)]{6} the expressions \eqref{56} for $C_1,C_2,D,G_1,G_2$
were also obtained, but starting from the second order $q$-difference operator
$L$ given in \eqref{55}, by which $K_1$ is represented in the polynomial
representation \cite[(17)]{6} of the Askey--Wilson algebra. Note that indeed
$h_n$ as given by \eqref{57} equals the eigenvalue of the eigenfunction $R_n(z)$ 
in \eqref{54}.

\subsection{A further symmetry of the Askey--Wilson algebra}
\label{60}
As observed in Section 2, nothing essentially changes for the
$q$-Verde-Star parameters
if $L$, $h_n$ and $g_n$ are simultaneously multiplied by a nonzero constant $\mu$,
or equivalently, if $a_1,a_2,g_1,g_2$
are multiplied by $\mu$. In the context of the $q$-Zhedanov algebra this means,
by \eqref{91},
that $D,G_1$ are multiplied by $\mu$ and $C_2,G_2,Q,\om$ are multiplied
by $\mu^2$.

Now, in the case of the Askey--Wilson algebra, take
$\mu=(q/(abcd))^\half=\td a^{-1}$ (with $\td a$ given by \eqref{53}.
Then it follows for the $q$-Zhedanov algebra generated by $K_1,K_2$ with
$K_1=\td a^{-1} L$ that, by \eqref{56} and \eqref{58}, $C_1,\td a^{-2} C_2,
\td a^{-1} D,\td a^{-1} G_1,
\td a^{-2} G_2$ and $\td a^{-2}\om$ are not only invariant under
permutations of $a,b,c,d$,
but also invariant if two parameters $e,f$ out of $a,b,c,d$ are sent
to $qe^{-1},qf^{-1}$, respectively (the symmetries already discussed
in Subsection \ref{59}).
Such symmetries for the Askey--Wilson DAHA were pointed out to me
by Marta Mazzocco. As was observed in Subsection \ref{59}, these symmetries
generate a group which is isomorphic to the Weyl group of the root system
$\textup{D}_4$. For some version of the Askey--Wilson algebra a symmetry
under the same group was
earlier observed in \cite[Section 4]{30}.

These further symmetries are relevant for the $q$-Zhedanov scheme
(Figure \ref{21} in the next section). If in \eqref{56}
some of the $C_1,C_2,D,G_1,G_2$ vanish under certain constraints on
$a,b,c,d$ then they will also vanish under transformed constraints
by application of a symmetry. For instance, the constraints $c=-a$, $d=-b$
give symmetric Askey--Wilson \eqref{28} and cause $D$ and $G_2$ to vanish.
Then also the constraints $c=-qa^{-1}$, $d=-qb^{-1}$ make $D$ and $G_2$ vanish.
Furthermore, these symmetries induce symmetries on limit cases of the
Askey--Wilson algebra.

The passage to the dual also becomes much nicer in this slightly adapted
Askey--Wilson algebra. Indeed note that by \eqref{57}, \eqref{61}, \eqref{56}
\begin{align*}
&b_1=a,\quad b_2=a^{-1},\quad \td a^{-1} a_1=\td a,\quad
\td a^{-1}a_2=\td a^{-1},\\
&x_k=aq^k+a^{-1}q^{-k},\quad \td a^{-1} h_k =\td aq^k+\td a^{-1}q^{-k},\\
&C_1=\td a^{-2} C_2=(q-q^{-1})^2,\quad
a^{-1}\td D=\td a^{-1} D,\quad
a^{-1}\td G_1=\td a^{-2} G_2,\quad
a^{-2}\td\om=\td a^{-2}\om,
\end{align*}
where $\td D$ is obtained from $D$ in \eqref{56} by replacing $a,b,c,d$ by
$\td a,\td b,\td c,\td d$, respectively, and similarly for $\td G_1$ and
$\td \om$. Thus $\td a^{-1} D,\td a^{-1} G_1,
\td a^{-2} G_2$ and $\td a^{-2}\om$ are also invariant under permutations
of $\td a,\td b,\td c,\td d$, which we already knew because of Proposition
\ref{62}.
\subsection{The dual algebras for continuous dual $q$-Hahn and big $q$-Jacobi}
\label{92}
For obtaining the continuous dual $q$-Hahn algebra put $d=0$
in \eqref{57}--\eqref{58}.
Then $e_4=0$ and $e_1,e_2,e_3$ can be considered as elementary symmetric
polynomials in $a,b,c$. We obtain
\begin{align}
&x_k=aq^k+a^{-1}q^{-k},\quad
h_k=q^{-k},\quad
g_k=q^{-2k+1}a^{-1}(1-abq^{k-1})(1-acq^{k-1})(1-q^k),\label{76}\sLP
&b_1=a,\quad
b_2=a^{-1}\quad
a_1=0,\quad
a_2=1,\quad
d_1=-q^{-1}abc,\quad
d_2=-(b+c+qa^{-1}),\label{77}\sLP
\begin{split}
&C_1 =(q-q^{-1})^2,\quad
C_2 =0,\quad
D =(1-q^{-1})^2(e_3+qe_1),\\
&G_1 =-(1-q^{-1})^2(1+q)(e_2+q),\quad
G_2 =-(1-q^{-1})^2(1+q)e_3.
\end{split}\label{78}
\end{align}
See also \cite[(50)]{6}.
Note that $U_n(z+z^{-1})=R_n(z;a,b,c,0\,|\,q)$
can be written in the form \eqref{51} with $x_k,h_k,g_k$ given by \eqref{76},
and that $h_n=q^{-n}$ is the eigenvalue of $L$ for $d=0$ in \eqref{54}.

Just as for the Askey--Wilson algebra nothing essentially changes if we multiply
$h_n,g_n$ by~$\mu$, and
$a_1,a_2,g_1,g_2$ by~$\mu$, and
$D,G_1$ by $\mu$, and $C_2,G_2$
by $\mu^2$ (here we have omitted $\om$).
If we now take $\mu=(q/(abc))^\half$ then the resulting
$C_1,\mu^2 C_2,
\mu D,\mu G_1,
\mu^2 G_2$ are not only invariant under permutations of $a,b,c$,
but also invariant if two parameters $e,f$ out of $a,b,c,$ are sent
to $qe^{-1},qf^{-1}$, respectively.

With $x_k,h_k,g_k$ given by \eqref{76} we get from \eqref{51} that
$\tilde U_m(y)=P_m(y;q^{-1}ab,ab^{-1},q^{-1}ac;q)$, a big $q$-Jacobi polynomial
as we saw it already in the duality \eqref{74}.

Below we give data associated with the big $q$-Jacobi polynomial
$U_n(x)=P_n(x;a,b,c;q)$.
\begin{align}
&x_k=q^{-k},\quad
h_k=q^{-k}+abq^{k+1},\quad
g_k=q^{1-2k}(1-aq^k)(1-cq^k)(1-q^k),\label{79}\\*[\smallskipamount]
&b_1=0,\quad
b_2=1,\quad
a_1=qab,\quad
a_2=1,\quad
d_1=-qac,\quad
d_2=-q(a+c+1),\label{80}\\*[\smallskipamount]
\begin{split}
&C_1 =0,\quad
C_2 =q^{-1}ab(1-q^2)^2,\quad
D =(1-q)^2(ab+ac+a+c),\\
&G_1 =-(1-q)^2(1+q)ac,\quad
G_2 =-(1-q)^2(1+q)a(ab+bc+b+c).
\end{split}\label{81}
\end{align}
If in \eqref{79}--\eqref{81} we multiply $h_k,g_k,a_1,a_2,d_1,d_2,D,G_1$ by
$(qab)^{-\half}$ and $C_2,G_2$ by $(qab)^{-1}$, and next replace
$a,b,c$ by $q^{-1}ab,ab^{-1},q^{-1}ac$ then we obtain the data which are dual to
the data \eqref{76}--\eqref{78}.

If in \eqref{79}--\eqref{81} we multiply $h_k,g_k,a_1,a_2,d_1,d_2,D,G_1$ by
$(qab)^{-\half}$ and $C_2,G_2$ by $(qab)^{-1}$, and if we multiply
$x_k,g_k,b_1,b_2,d_1,d_2,D,G_2$ by $(ab^{-1}c^2)^{-\frac14}$ and $C_1,G_1$ by
$(ab^{-1}c^2)^{-\half}$ then
$U_n(x)=P_n\big((ab^{-1}c^2)^{\frac14}x;a,b,c;q\big)$. Moreover
the big $q$-Jacobi algebra then has the invariances corresponding to the ones we
observed for the continuous dual Hahn algebra. We see that
$C_1,C_2,D,G_1,G_2$ are invariant under the transformations
$(a,b,c)\to(c,abc^{-1},a)$, $(a,b,c)\to (abc^{-1},c,b)$ and
$(a,b,c)\to (a,b^{-1},b^{-1}c)$. This nicely corresponds with
equalities \eqref{67}, \eqref{68} and with the symmetry observed in connection
with \eqref{74}.
\section{The $q$-Zhedanov scheme}
\label{89}
It will turn out that vanishing of some of the coefficients $C_1,C_2,D,G_1,G_2$
in the $q$-Zhedanov algebra with  relations \eqref{12} is a
characterizing property
of corresponding polynomials in the $q$-Askey scheme. An arrow then 
will indicate that one more of the five constants becomes zero.
We arrange the five coefficients in an array
\begin{equation}\label{20}
\begin{matrix}C_1&&C_2\\&D&\\G_1&&G_2\end{matrix}
\end{equation}
By \eqref{19} we can pass to the dual $q$-Zhedanov algebra by reflection of the
array with respect to the vertical axis.

We will replace coefficients in this array by symbols \black\ or \white, where
\black\ denotes any value (including zero) and \white\ denotes zero.
From the explicit data for the $q$-Verde-Star parameters of the
families in the $q$-Askey scheme as given in \cite[Appendix A]{1} we can obtain
the corresponding vanishing pattern for the array \eqref{20} by means of
\eqref{14}. Then connect the
arrays by arrows such that in the direction of the arrow at least one \black\
is turned into \white. The resulting graph is in Figure~\ref{21}.
\setlength{\unitlength}{2.35mm}
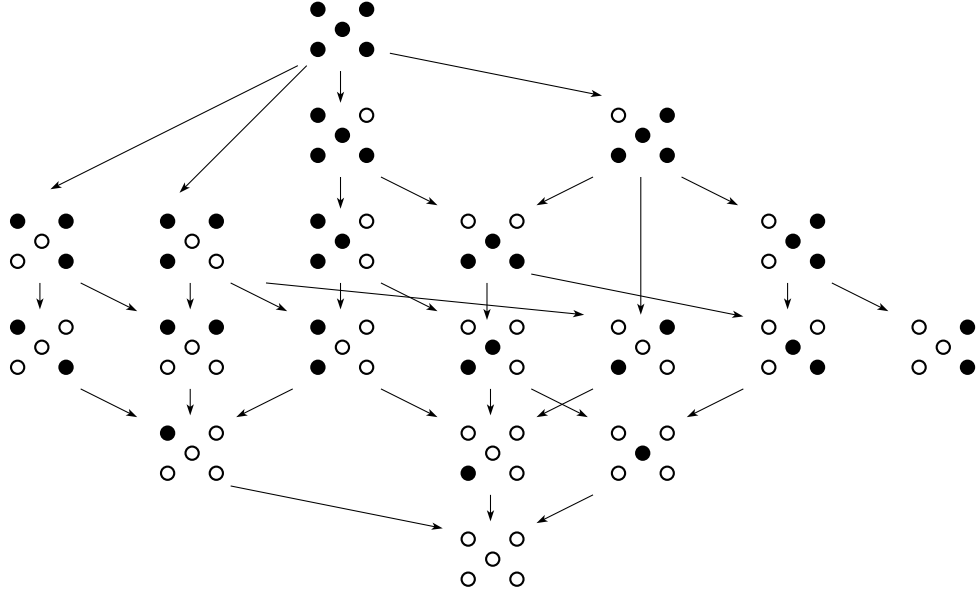
\begin{figure}[ht]
\centering
\hskip-3cm
\begin{picture}(30,40)
\put(11.5,37)
{\text{$\begingroup\setlength\arraycolsep{1.7pt}
\begin{matrix}\black&&\black\\[-6pt]
&\black&\\[-6pt]\black&&\black\end{matrix}
\endgroup$}}
\put(10.8,35.3) {\vector(-2,-1){14}}
\put(11.3,35.3) {\vector(-1,-1){7.2}}
\put(13.2,35) {\vector(0,-1){1.8}}
\put(16,36) {\vector(5,-1){12}}
\put(11.5,31)
{\text{$\begingroup\setlength\arraycolsep{1.7pt}
\begin{matrix}\black&&\white\\[-6pt]
&\black&\\[-6pt]\black&&\black\end{matrix}
\endgroup$}}
\put(13.2,29) {\vector(0,-1){1.8}}
\put(15.5,29) {\vector(2,-1){3.2}}
\put(28.5,31)
{\text{$\begingroup\setlength\arraycolsep{1.7pt}
\begin{matrix}\white&&\black\\[-6pt]
&\black&\\[-6pt]\black&&\black\end{matrix}
\endgroup$}}
\put(27.5,29) {\vector(-2,-1){3.2}}
\put(30.2,29) {\vector(0,-1){7.8}}
\put(32.5,29) {\vector(2,-1){3.2}}
\put(-5.5,25)
{\text{$\begingroup\setlength\arraycolsep{1.7pt}
\begin{matrix}\black&&\black\\[-6pt]
&\white&\\[-6pt]\white&&\black\end{matrix}
\endgroup$}}
\put(-3.8,23) {\vector(0,-1){1.5}}
\put(-1.5,23) {\vector(2,-1){3.2}}
\put(3,25)
{\text{$\begingroup\setlength\arraycolsep{1.7pt}
\begin{matrix}\black&&\black\\[-6pt]
&\white&\\[-6pt]\black&&\white\end{matrix}
\endgroup$}}
\put(4.7,23) {\vector(0,-1){1.5}}
\put(7,23) {\vector(2,-1){3.2}}
\put(9,23) {\vector(10,-1){18}}
\put(11.5,25)
{\text{$\begingroup\setlength\arraycolsep{1.7pt}
\begin{matrix}\black&&\white\\[-6pt]
&\black&\\[-6pt]\black&&\white\end{matrix}
\endgroup$}}
\put(13.2,23) {\vector(0,-1){1.5}}
\put(15.5,23) {\vector(2,-1){3.2}}
\put(20,25)
{\text{$\begingroup\setlength\arraycolsep{1.7pt}
\begin{matrix}\white&&\white\\[-6pt]
&\black&\\[-6pt]\black&&\black\end{matrix}
\endgroup$}}
\put(21.5,23) {\vector(0,-1){2.1}}
\put(24,23.5) {\vector(5,-1){12}}
\put(37,25)
{\text{$\begingroup\setlength\arraycolsep{1.7pt}
\begin{matrix}\white&&\black\\[-6pt]
&\black&\\[-6pt]\white&&\black\end{matrix}
\endgroup$}}
\put(38.5,23) {\vector(0,-1){1.5}}
\put(41,23) {\vector(2,-1){3.2}}
\put(-5.5,19)
{\text{$\begingroup\setlength\arraycolsep{1.7pt}
\begin{matrix}\black&&\white\\[-6pt]
&\white&\\[-6pt]\white&&\black\end{matrix}
\endgroup$}}
\put(-1.5,17) {\vector(2,-1){3.2}}
\put(3,19)
{\text{$\begingroup\setlength\arraycolsep{1.7pt}
\begin{matrix}\black&&\black\\[-6pt]
&\white&\\[-6pt]\white&&\white\end{matrix}
\endgroup$}}
\put(4.7,17) {\vector(0,-1){1.5}}
\put(11.5,19)
{\text{$\begingroup\setlength\arraycolsep{1.7pt}
\begin{matrix}\black&&\white\\[-6pt]
&\white&\\[-6pt]\black&&\white\end{matrix}
\endgroup$}}
\put(10.5,17) {\vector(-2,-1){3.2}}
\put(15.5,17) {\vector(2,-1){3.2}}
\put(20,19)
{\text{$\begingroup\setlength\arraycolsep{1.7pt}
\begin{matrix}\white&&\white\\[-6pt]
&\black&\\[-6pt]\black&&\white\end{matrix}
\endgroup$}}
\put(21.7,17) {\vector(0,-1){1.5}}
\put(24,17) {\vector(2,-1){3.2}}
\put(28.5,19)
{\text{$\begingroup\setlength\arraycolsep{1.7pt}
\begin{matrix}\white&&\black\\[-6pt]
&\white&\\[-6pt]\black&&\white\end{matrix}
\endgroup$}}
\put(27.5,17) {\vector(-2,-1){3.2}}
\put(37,19)
{\text{$\begingroup\setlength\arraycolsep{1.7pt}
\begin{matrix}\white&&\white\\[-6pt]
&\black&\\[-6pt]\white&&\black\end{matrix}
\endgroup$}}
\put(36,17) {\vector(-2,-1){3.2}}
\put(45.5,19)
{\text{$\begingroup\setlength\arraycolsep{1.7pt}
\begin{matrix}\white&&\black\\[-6pt]
&\white&\\[-6pt]\white&&\black\end{matrix}
\endgroup$}}
\put(3,13)
{\text{$\begingroup\setlength\arraycolsep{1.7pt}
\begin{matrix}\black&&\white\\[-6pt]
&\white&\\[-6pt]\white&&\white\end{matrix}
\endgroup$}}
\put(7,11.5) {\vector(5,-1){12}}
\put(20,13)
{\text{$\begingroup\setlength\arraycolsep{1.7pt}
\begin{matrix}\white&&\white\\[-6pt]
&\white&\\[-6pt]\black&&\white\end{matrix}
\endgroup$}}
\put(21.7,11) {\vector(0,-1){1.5}}
\put(28.5,13)
{\text{$\begingroup\setlength\arraycolsep{1.7pt}
\begin{matrix}\white&&\white\\[-6pt]
&\black&\\[-6pt]\white&&\white\end{matrix}
\endgroup$}}
\put(27.5,11) {\vector(-2,-1){3.2}}
\put(20,7)
{\text{$\begingroup\setlength\arraycolsep{1.7pt}
\begin{matrix}\white&&\white\\[-6pt]
&\white&\\[-6pt]\white&&\white\end{matrix}
\endgroup$}}
\end{picture}
\vspace*{-1.5cm}
\caption{The $q$-Zhedanov scheme. In referring to this
rows are numbered top to bottom, and in each
row the elements are labeled by the letters a,b,c,$\ldots$ from left to right.}
\label{21}
\end{figure}

Let us number the rows in the scheme from top to bottom by 1 to 6.
In each row list the successive arrays from left to right by a,~b,~$\ldots$\;.

Note the following dualities in Figure \ref{21}:\\
\textbf{1a}, \textbf{3d}, \textbf{4b}, \textbf{5c}, \textbf{6a}
are self-dual.\\
$\textbf{2a}\leftrightarrow\textbf{2b}$,
$\textbf{3a}\leftrightarrow\textbf{3b}$,
$\textbf{3c}\leftrightarrow\textbf{3e}$,
$\textbf{4a}\leftrightarrow\textbf{4e}$,
$\textbf{4c}\leftrightarrow\textbf{4g}$,
$\textbf{4d}\leftrightarrow\textbf{4f}$
are dual pairs.\\
The duals of \textbf{5a} and \textbf{5b} are not included in the scheme,
but will be discussed in Section \ref{88}.

The arrays in Figure \ref{21} correspond with families in 
the $q$-Askey scheme as given in the list below.
Usually we mention only one family for an array, but also subfamilies still
corresponding to that array. There may be more families
for the same array, by $q\leftrightarrow q^{-1}$ exchange, and there may be
discrete families. If one wants to add these, one can use \cite[Ch.~14]{3}.
We do not bother
at all about orthogonality properties. The only things that matter
are that the polynomials are eigenfunctions of a second order $q$-difference
operator and satisfy a three-term recurrence relation.

As much as possible formulas for the occurring families in the list
will be given in terms of Askey--Wilson $R_n(z;a,b,c,d\,|\,q)$,
$x=z+z^{-1}$
(see \eqref{27}) or big $q$-Jacobi $P_n(x;a,b,c;q)$ (see \eqref{48}).
Otherwise a formula is given as a limit case of $R_n$ or $P_n$ or
as some ${}_r\phi_s$ $q$-hypergeometric
function. Everything is up to constant factors.
\mLP
\textbf{1a.} Askey--Wilson $R_n(z;a,b,c,d\,|\,q)$, see \eqref{27}.
\mLP
\textbf{2a.} continuous dual $q$-Hahn $R_n(z;a,b,c,0\,|\,q)$, see \eqref{83}.
\mLP
\textbf{2b.} big $q$-Jacobi $P_n(x;a,b,c;q)$, see \eqref{48}.
\mLP
\textbf{3a.} continuous $q^2$-Jacobi
$R_n(z;a,b,q^\half,-q^\half\,|\,q)$, see \eqref{33}.
\mLP
\textbf{3b.} symmetric Askey--Wilson $R_n(z;a,b,-a,-b\,|\,q)$
(see \eqref{28}) or $R_n(z;a,b,-qa^{-1},-qb^{-1}\,|\,q)$.
\mLP
\textbf{3c.} Al-Salam--Chihara $R_n(z;a,b,0,0\,|\,q)$ (see \eqref{82})
with subfamily continuous big $q$-Hermite $R_n(z;a,0,0,0\,|\,q)$, see
\cite[(14.8.1)]{3}.
\mLP
\textbf{3d.} big $q$-Laguerre $P_n(x;a,0,c;q)$, see \eqref{84}.
\mLP
\textbf{3e.} little $q$-Jacobi $P_n(x;a,b,0;q)$ (see \eqref{70})
with limit family $q$-Bessel\\
$\lim_{b\to0} P_n(-ab^{-1}x;-q^{-1}ab^{-1},b,0;q)=
{}_2\phi_1(q^{-n},-q^na;0;q,qx)$ (see \cite[(14.22.1)]{3}).
For $a=q^{-N-1}$, $b=-p$ little $q$-Jacobi becomes
$q$-Krawtchouk $P_n(x;q^{-N-1},-p,0;q)$ (see \cite[(14.15.1)]{3}).
\mLP
\textbf{4a.} continuous $q^2$-Laguerre
$R_n(z;q^\half,-q^\half,a,0\,|\,q)$, see \eqref{34}.
\mLP
\textbf{4b.} continuous $q^2$-ultraspherical
$R_n(z;q^\half,-q^\half,a,-a\,|\,q)$ (see \eqref{35})\\
or $R_n(z;q^\half,-q^\half,a,-qa^{-1}\,|\,q)$.
\mLP
\textbf{4c.} symmetric Al-Salam--Chihara $R_n(z;a,-a,0,0\,|\,q)$
(see \eqref{85}) with limit family continuous $q$-Hermite
$\lim_{a\to0}a^{-n} R_n(z;a,-a,0,0\,|\,q)=H_n\big(\thalf(z+z^{-1})\,|\,q\big)$
(see \eqref{32}).
\mLP
\textbf{4d.} Al-Salam--Carlitz I
$\lim_{a\to0}a^{-n}P_n(qax;a,0,ac;q)=
(-c)^n q^{\half n(n+1)} {}_2\phi_1(q^{-n},x^{-1};0;q,qc^{-1}x)$
(see \cite[(14.24.1)]{3}).
\mLP
\textbf{4e.} symmetric big $q$-Jacobi $P_n(x;a,a,-a;q)$ (see \eqref{86})
or $P_n(x;-1,-1,c;q)$.
\mLP
\textbf{4f.} $q$-Laguerre
$\lim_{a\to-\iy} P_n(-x;a,b,0;q)=(qb;q)_n
{}_1\phi_1(q^{-n};qb;q,-q^{n+1}bx)$ (see \eqref{71}, \eqref{87},
\cite[(14.21.1), (14.12.13)]{3})
with
limit family Stieltjes--Wigert
$\lim_{b\to0}\lim_{a\to-\iy} P_n(-b^{-1}x;a,b,0;q)=
{}_1\phi_1(q^{-n};0;q,-q^{n+1}x)$ (see \cite[(14.27.1)]{3}).
$q$-Charlier is $q$-Laguerre with $x,b$ replaced by $-x,-b^{-1}$
(see \cite[(14.23.1)]{3}).
This case also contains  little $q$-Laguerre, which can be essentially identified
with $q^{-1}$-Laguerre, see \cite[p.521]{3}.
\mLP
\textbf{4g.} special little $q$-Jacobi $P_n(a,-1,0;q)$. For $a=q^{-N-1}$ this
is $q$-Krawtchouk with $p=q^{-N}$, see \cite[(14.15.1)]{3}.
\mLP
\textbf{5a.} continuous $q^2$-Hermite $R_n(z;q^\half,-q^\half,0,0\,|\,q)$,
see \eqref{36}.
\mLP
\textbf{5b.} discrete $q$-Hermite I
$\lim_{a\to0}a^{-n}P_n(qax;a,0,-a;q)=
q^{\half n(n+1)} {}_2\phi_1(q^{-n},x^{-1};0;q,-qx)$
(see \cite[(14.28.1)]{3}).
\mLP
\textbf{5c.} $x^n(x^{-1};q)_n$
\mLP
\textbf{6a.} $x^n$
\paragraph{Remarks}\quad\sLP
\textbf{1.} For obtaining this list we started with the $q$-Verde-Star data
given in \cite[Appendix A]{1}. From these data, with $a_0$ and $b_0$ being
put to zero, the $q$-Zhedanov coefficients $C_1,C_2,D,G_1,G_2$ can be
computed by \eqref{14},
and it can be read off which of these five coefficients vanish.
Moreover, for a given family, it can be seen for which constraints on
the parameters a further coefficient will vanish.
\mLP
\textbf{2.} It is interesting to compare the scheme in Figure \ref{21} with
the $q$-Verde-Star scheme \cite[Figure~1]{1}:
\begin{itemize}
\item
First, the present scheme is more
compact because it does not discern between $q$ and $q^{-1}$, and 
also because it does not
bother about the Newton polynomials $v_k$ in which $u_n$ is expanded,
so each family occurs in at most one place in the scheme.
\item
Furthermore, some separate families in \cite[Figure~1]{1} here merge as
subfamilies or limit families with other families (continuous big $q$-Hermite,
$q$-Bessel and Stieltjes--Wigert).
\item
On the other hand, the present scheme is richer because many subfamilies
here get an independent status which they did not have in the
$q$-Verde-Star scheme. We mention symmetric Askey--Wilson,
symmetric Al-Salam--Chihara, symmetric big $q$-Jacobi,
continuous $q^2$-Jacobi, continuous $q^2$-Laguerre,
continuous $q^2$-ultraspherical, discrete $q$-Hermite and
$q$-Krawtchouk with $p=q^{-N}$.
\item
Continuous $q$-Hermite could not be included in the $q$-Verde-Star scheme,
because it cannot be written in the form \eqref{1}, but a
$q$-Zhedanov algebra can be associated with it.
Therefore, in the present scheme it occurs, but only as a limit family of
symmetric Al-Salam--Chihara.
However, continuous $q^2$-Hermite, obtained by the quadratic transformation
\eqref{36}, has its own place in our scheme. We also see from \eqref{36}
that continuous $q^2$-Hermite will live in the
$q$-Verde-Star scheme as a subfamily of Askey--Wilson, but not with a 
separate place in that scheme. To some extent this overcomes the defect that
continuous $q$-Hermite could not be handled by Verde-Star \cite{9}.
\end{itemize}
\textbf{3.}
The cases \textbf{5c} and \textbf{6} are degenerate. They are eigenfunctions
of a first order $q$-difference operator and they satisfy a two-term
recurrence relation.
\mLP
\textbf{4.}
It turns out that with each of the 32 possible arrays a family of polynomials
can be associated, but we have not included all of these in the scheme,
because some families did not seem to be of sufficient interest.
The following cases were missing in the scheme.
\sLP
$\begingroup\setlength\arraycolsep{1.7pt}
\begin{matrix}\black&&\black\\[-6pt]
&\white&\\[-6pt]\black&&\black\end{matrix}
\endgroup$\quad
Askey--Wilson $R_n(z;a,b,c,d\,|\,q)$ with
$a b c + a b d + a c d + b c d + q (a + b + c + d)=0$.
\sLP
$\begingroup\setlength\arraycolsep{1.7pt}
\begin{matrix}\black&&\black\\[-6pt]
&\black&\\[-6pt]\white&&\black\end{matrix}
\endgroup$\quad
idem with $a b c d + q (a b + a c + a d + b c + b d + c d) + q^2=0$.
\sLP
$\begingroup\setlength\arraycolsep{1.7pt}
\begin{matrix}\black&&\black\\[-6pt]
&\black&\\[-6pt]\black&&\white\end{matrix}
\endgroup$\quad
idem with $a b c d (a + b + c + d) + q (a b c + a b d + a c d + b c d)=0$.
\sLP
$\begingroup\setlength\arraycolsep{1.7pt}
\begin{matrix}\black&&\black\\[-6pt]
&\black&\\[-6pt]\white&&\white\end{matrix}
\endgroup$\quad
idem with the previous two constraints.
\sLP
$\begingroup\setlength\arraycolsep{1.7pt}
\begin{matrix}\black&&\white\\[-6pt]
&\white&\\[-6pt]\black&&\black\end{matrix}
\endgroup$\quad
continuous dual $q$-Hahn $R_n(z;a,b,c,0\,|\,q)$ with
$abc+q(a+b+c)=0$.
\sLP
$\begingroup\setlength\arraycolsep{1.7pt}
\begin{matrix}\black&&\white\\[-6pt]
&\black&\\[-6pt]\white&&\black\end{matrix}
\endgroup$\quad
idem with $ab+ac+bc+q=0$.
\sLP
$\begingroup\setlength\arraycolsep{1.7pt}
\begin{matrix}\white&&\black\\[-6pt]
&\white&\\[-6pt]\black&&\black\end{matrix}
\endgroup$\quad
big $q$-Jacobi $P_n(x;a,b,c;q)$ with $a+c+a(b+c)=0$.
\sLP
$\begingroup\setlength\arraycolsep{1.7pt}
\begin{matrix}\white&&\black\\[-6pt]
&\black&\\[-6pt]\black&&\white\end{matrix}
\endgroup$\quad
idem with $b+c+b(a+c)=0$.
\sLP
$\begingroup\setlength\arraycolsep{1.7pt}
\begin{matrix}\black&&\white\\[-6pt]
&\black&\\[-6pt]\white&&\white\end{matrix}
\endgroup$\quad
Al-Salam--Chihara $R_n(z;a,b,0,0\,|\,q)$ with $ab+q=0$.
\sLP
$\begingroup\setlength\arraycolsep{1.7pt}
\begin{matrix}\white&&\black\\[-6pt]
&\black&\\[-6pt]\white&&\white\end{matrix}
\endgroup$\quad
special little $q$-Jacobi $P_n(x;-1,b,0;q)$.
\sLP
$\begingroup\setlength\arraycolsep{1.7pt}
\begin{matrix}\white&&\white\\[-6pt]
&\white&\\[-6pt]\black&&\black\end{matrix}
\endgroup$\quad
big $q$-Laguerre $P_n(x;a,0,c;q)$ with $a+c+ac=0$,
\sLP
$\begingroup\setlength\arraycolsep{1.7pt}
\begin{matrix}\white&&\black\\[-6pt]
&\white&\\[-6pt]\white&&\white\end{matrix}
\endgroup$\quad
dual continuous $q^2$-Hermite $\dstyle\qhyp32{q^{-n},q^{n+1},x}{-q,0}{q,q}$
(see Section \ref{88}).
\sLP
$\begingroup\setlength\arraycolsep{1.7pt}
\begin{matrix}\white&&\white\\[-6pt]
&\white&\\[-6pt]\white&&\black\end{matrix}
\endgroup$\quad
special $q$-Charlier $\dstyle\qhyp21{q^{-n},-x}0{q,-q^{n+1}}$ 
(see Section \ref{88}).
\mLP
\textbf{5.}
In Subsections \ref{60}, \ref{92} we already observed symmetries for
$q$-Zhedanov algebras which are not satisfied by the corresponding
Askey--Wilson, continuous dual $q$-Hahn and big $q$-Jacobi polynomials. 
Such symmetries will also be present in algebras associated with families
in the third and fourth row of Figure 1.
\section{Duals of continuous $q^2$-Hermite and discrete $q$-Hermite I}
\label{88}
In this section we will phrase the various dualities in the notation of
\eqref{50}--\eqref{52}.
\paragraph{Continuous $q^2$-Laguerre $\longleftrightarrow$ special
big $q$-Jacobi}\quad\\
For obtaining the dual of continuous $q$-Hermite we first consider a suitable
dual of continuous $q^2$-Laguerre and then specialize the parameter.

Put $x_k=q^{k+\half}+q^{-k-\half}$, $h_k=q^{-k}$,
$g_k=(q^\half-aq^k)(q^{-2k}-1)$. 
Then by \eqref{14} we have
$C_1=(q-q^{-1})^2$, $C_2=D=G_1=0$, $G_2=(1+q^{-1})(1-q)^2a$.
Also \eqref{50}--\eqref{52} take the form
\begin{equation}
\begin{split}
U_n(z+z^{-1})=\qhyp32{q^{-n},q^\half z,q^\half z^{-1}}{q^\half a,-q}{q,q},
\qquad\tilde U_m(y)=\qhyp32{y,q^{-m},q^{m+1}}{q^\half a, -q}{q,q},\\
U_n(q^{m+\half}+q^{-m-\half})=\tilde U_m(q^{-n}).
\end{split}
\label{63}
\end{equation}
By \eqref{27}, \eqref{34} and \eqref{48} we can identify the dual polynomials
$U_n(z+z^{-1})$ and $\tilde U_m(y)$ as continuous $q^2$-Laguerre and
special big $q$-Jacobi:
\begin{equation}
U_n(z+z^{-1})=R_n(z;q^\half,-q^\half,a,0\,|\,q),\quad
\tilde U_m(y)=P_m(y;-1,-1,q^{-\half} a;q).
\label{64}
\end{equation}

It is precisely this dual pair which allows the specialization $a=0$ which
we need. The $\tilde U_m$ belong to the same array
\textbf{4e} in Figure 1
as symmetric big $q$-Jacobi. In fact, with the choice
$x_k=aq^k+a^{-1}q^{-k}$, $h_k=q^{-k}$,
$g_k=q^{1-2k}a^{-1}(1-q^k)(1-a^2 q^{2k-1})$ we would have obtained the
same values for $C_1,C_2,D,G_1,G_2$, but
$U_n(z+z^{-1})=R_n(z;a,q^\half,-q^\half,0\,|\,q)$ and for
$\tilde U_m(y)$ the symmetric big $q$-Jacobi polynomial
$P_m(y;q^{-\half}a,q^{-\half}a,-q^{-\half}a;q)$.

It is also remarkable that for both specializations of big $q$-Jacobi
belonging to array \textbf{4e} the specialization of the three-term
recurrence relation \eqref{69} gives (for $n=1,2,\ldots$)
$xP_n(x)$ as a linear combination of $P_{n+1}(x)$ and $P_{n-1}(x)$ without
a term for $P_n(x)$. Still only for symmetric big $q$-Jacobi
$P_n(x;a,a,-a;q)$ we have $P_n(-x)=(-1)^n P_n(x)$, but for
$P_n(x;-1,-1,c;q)$ the starting value $P_1(x)$ already lacks this symmetry.
\paragraph{Continuous $q^2$-Hermite $\longleftrightarrow$ very special
big $q$-Jacobi}\quad\\
Just put $a=0$ in the above case.

Note that the ``dual continuous $q^2$-Hermite polynomials''
$P_n(x)=P_n(x;-1,-1,0;q)$ cannnot be orthogonal with respect to a positive
orthogonality measure. This is seen from the signs of the coefficients
in the recurrence relation (a specialization of \eqref{69})
\begin{equation}
xP_n(x)=\frac1{1-q^{2n+1}}\,P_{n+1}(x)-\frac{q^{2n+1}}{1-q^{2n+1}}\,
P_{n-1}(x),\quad n=1,2,\ldots\;.
\end{equation}
\paragraph{Al-Salam--Carlitz I $\longleftrightarrow$ $q$-Charlier}\quad\\
For obtaining the dual of discrete $q$-Hermite I we first consider the
dual families of Al-Salam--Carlitz I and $q$-Charlier and then specialize
the parameter.

Put $x_k=q^k$, $h_k=q^{-k}$, $g_k=a(q^{-k}-1)$. 
Then $C_1=C_2=0$, $D=q^{-1}(1-q)^2(1-a)$,
$G_1=(1-q^{-1})^2(1+q)a$, $G_2=0$.
Also \eqref{50}--\eqref{52}
take the form
\begin{equation}
\begin{split}
U_n(x)=\qhyp21{q^{-n},x^{-1}}0{q,-qa^{-1}x},\qquad
\wt U_m(y)=\qhyp21{y,q^{-m}}0{q,-q^{m+1}a^{-1}},\\
U_n(q^m)=\wt U_m(q^{-n}).
\end{split}
\end{equation}
Furthermore,
by \cite[(14.24.1), (14.23.1)]{3}, the expressions in terms of
Al-Salam--Carlitz I and $q$-Charlier are:
\begin{equation}
U_n(x)=a^n q^{-\half n(n-1)} U_n^{(-a)}(x;q),\qquad
\wt U_m(y)=C_m(y;a;q).
\end{equation}
\paragraph{Discrete $q$-Hermite I $\longleftrightarrow$ special
$q$-Charlier}\quad\\
Just put $a=1$ in the above case. Then
$D$ will also vanish and,
by \cite[(14.28.1), (14.23.1)]{3}, the expressions in terms of
discrete $q$-Hermite I and $q$-Charlier are:
\begin{equation}
U_n(x)=q^{-\half n(n-1)} h_n(x;q),\qquad
\wt U_m(y)=C_m(y;1;q)..
\end{equation}
\section{Charting the $q$-Askey scheme: conclusions and questions}
\label{95}
From Sections \ref{90} and \ref{89} we can conclude that a system $\{p_n\}$,
belonging to a family in the $q$-Askey scheme for some values of the
parameters for that family, and considered up to dilation of the
independent variable, gives rise to values of
$C_1,C_2,D,G_1,G_2,\om$ up to multiplication by certain powers of nonzero
constants $\mu$ and $\rho$ as in \eqref{91}. In view of
Subsections \ref{60}, \ref{92} and Remark~5 in Section \ref{89} the
system $\{p_n\}$
will not be uniquely determined by the equivalence class of
$C_1,C_2,D,G_1,G_2,\om$. If the $p_n$ are Askey--Wilson polynomials then we
may conjecture from Subsection \ref{60} that its parameters
$\pm(a,b,c,d)$ are determined
by the equivalence class of $C_1,C_2,D,G_1,G_2,\om$ up to permutations
and up to sending any $e,f$ from $a,b,c,d$ to $q/e$, $q/f$, but even that
is not yet clear.

Still it is of interest to let the equivalence classes \eqref{91} form a
complex manifold. Here we impose that $\om$ and at least one of the other five
coordinates are nonzero. Then we can put $\om$ and another nonvanishing
coordinate equal to 1. This will fix $\mu$ and $\rho$, possibly up to sign.
Then use the other four coordinates as local coordinates. There will still
be point pairs which have to be identified.

Another question, more a matter of taste and of convenience than of
mathematical rigor, is whether the usual version \cite[p.414]{3} of the
$q$-Askey scheme should be adapted in view of the $q$-Verde-Star scheme
\cite{1} and the $q$-Zhedanov scheme Figure 1. I would say that an adapted 
scheme should not distinguish between $q$ and $q^{-1}$
variants of the same family. Neither should it distinguish between two cases
of the same family with expansions in terms of different Newton polynomials.
Interesting as these features of the $q$-Verde-Star scheme may be,
they make the graphical display too large and too detailed.
For the rest the adapted scheme should be the union
of the $q$-Verde-Star and the $q$-Zhedanov scheme, since each of the schemes
highlights some families which deserve special attention but are not
in the other scheme. It would be convenient to have one scheme which starts
with Askey--Wilson and another scheme which starts with $q$-Racah. In the lower
rows these schemes would have several families in common.

Transcendental functions might also be associated with the $q$-Zhedanov scheme,
as we will argue now.
The $q$-Zhedanov algebra yielding the coefficients $C_1,C_2,D,G_1,G_2$ by
\eqref{12} is generated by a second order $q$-difference operator $K_1$
and the operator $K_2$ which is multiplication by~$x$. Polynomials
associated with the algebra occur as eigenfunctions of $K_1$ for special
eigenvalues. However, the algebra and its coefficients are independent
of the choice of the eigenvalues. Other eigenfunctions of $K_1$ on another
domain may yield generalized orthogonal systems leading to integral
transforms. Notably, in addition to Askey--Wilson polynomials, there are
Askey--Wilson functions, see Koelink and Stokman \cite{12}. The same authors
have in \cite[Figure 1.2]{11} a scheme, somewhat analogous to the $q$-Askey
scheme, with limit cases of the Askey--Wilson functions. It would be 
interesting to match our scheme with a possibly extended version of the
scheme in~\cite{11}.
Passage to the transcendental case will also be needed for realizing some symmetries
of the $q$-Zhedanov algebra. For instance, the symmetry
$(a,b,c,d)\to(a,b,qd^{-1},qc^{-1})$ of the Askey--Wilson algebra (see Subsection
\ref{60}) is, by \cite[(3.2)]{11}, visible for Askey--Wilson functions, while this
symmetry is not present if we restrict ourselves to the polynomial case.

Schemes for nonsymmetric Askey--Wilson polynomials and their subfamilies, and
for the corresponding DAHAs might be another follow-up topic. A Verde-Star type
approach for nonsymmetric polynomials may be promising. However, in view of the
cases already considered in~\cite{6},
there does  not seem to be a uniform approach for the DAHA
as we do have by \eqref{12} for the $q$-Zhedanov algebras. 
\paragraph{Acknowledgements}
I thank Marta Mazzocco for communicating to me about the symmetry
$(a,b,c,d)\to(a,b,qd^{-1},qc^{-1})$ of the Askey--Wilson algebra.
I also thank Raimundas Vidunas for pointing me to \cite{19} and \cite{21}, and
an anonymous referee for helpful comments.

\quad\\
\begin{footnotesize}
\begin{quote}
{ T. H. Koornwinder, Korteweg-de Vries Institute, University of
 Amsterdam,\\
 P.O.\ Box 94248, 1090 GE Amsterdam, The Netherlands;

\vspace{\smallskipamount}
email: }\url{thkmath@xs4all.nl}
\end{quote}
\end{footnotesize}

\end{document}